\documentstyle{amsppt}
\document
\NoRunningHeads
\NoBlackBoxes
\TagsOnRight
\mag=1000
\hoffset=.15in
\pageno=1
%\input psbox.tex 
%\psfortextures
\define\subord{\hbox{$\triangleright\kern-.17cm\blacktriangleleft$}}
\define\subordd{\hbox{$\blacktriangleright\kern-.17cm\triangleleft$}}    
\catcode`\@=11
\redefine\logo@{}
\catcode`\@=\active

\hfill {\bf Swan 00-8}

\heading{\bf{A GENERALISED HOPF ALGEBRA FOR SOLITONS}}
\endheading

\heading{Falleh R. Al-Solamy\\
Department of Mathematics,\\ Faculty of Science, \\King Abdul Aziz
University, \\P. O. Box 80015,
\\Jeddah 21589, Saudi Arabia.\\ \& \\Edwin J. Beggs\\
Department of Mathematics,\\ University of Wales\\ Swansea SA2 8PP, Wales, UK.}
\endheading

\heading{Abstract}\endheading
This paper considers a generalisation of the idea of a Hopf algebra in which a commutative ring
replaces the field in the unit and counit. It is motivated by an example from the inverse scattering
formalism for solitons. 
We begin with the corresponding idea for groups, where the concept of the identity is altered. 

\heading{1. Introduction}\endheading
Group factorisation plays a vital part in the inverse scattering procedure [2,10].
For example the Riemann-Hilbert problem is a (not quite exact)
factorisation of group valued functions on the real line into functions
analytic on the lower half plane times
functions analytic on the upper half plane.
However there is a problem, a group valued function which is analytic on the
lower half plane need not have an inverse which is analytic there.
On the Lie algebra level all is well since any smooth loop which is
uniformly sufficiently close to the identity  and is analytic on the
lower half plane has an inverse which is also analytic on the
lower half plane. To avoid the problem, we look at the Lie algebras
or a neighbourhood of a group near the identity. This corresponds
in inverse scattering to looking at solutions not too far from the vacuum.
However the soliton solutions
for many integrable systems are characterised by meromorphic
loops, and there the factors are very definitely not closed under inverse.
 For example, if we take the meromorphic function given by
(for $P\in M_n(\Bbb C)$ a Hermitian projection matrix and $P^{\bot}=1-P$)
$$
\phi(\lambda)= P^{\bot} + \frac{\lambda - \overline\alpha}{\lambda - \alpha} P\ ,\tag1.1
$$ which has pole at $\alpha$ in the upper half plane, 
its inverse is given by 
$$\phi^{-1}(\lambda)= P^{\bot} + \frac{\lambda - \alpha}{\lambda - \overline\alpha} P\ ,\tag1.2
$$ which has pole at $\overline\alpha$ in the lower half plane.

The meromorphic loops which specify the solitons in a classical integrable system
are not uniquely defined, there are `vacuum loops' which can be added
without generating any extra solitons. However these can be thought of as allowing
soliton-antisoliton
pair creation in the integrable system. When the system is quantised
a vacuum loop could be perturbed into a soliton-antisoliton pair by a slight
movement of the pole positions. Effectively the vacuum loops get round the
problem which would occur if we could count solitons by the number of poles.
The number of poles cannot be changed by a small perturbation, so without the
vacuum loops soliton-antisoliton pair production might seem impossible. 
To calculate the total quantum energy and momentum
for a soliton the contributions for these vacuum loops would have to be 
added. 
This quantum correction is observed in the Sine-Gordon model, calculated by other methods [13].

The existence of the vacuum loops and the fact that the upper and lower factors
for the meromorphic loop group are not groups are related, and are both taken into
account in the ideas described in this paper of an `almost group' and 
matched pairs of almost groups. This naturally leads on to the idea of an `almost 
Hopf algebra', in which the unit and counit map are modified to use a commutative
algebra instead of the ground field. From the discussion above, the commutative algebra 
would arise from the vacuum loops. A group factorisation into a subgroup
(a group doublecross product) is well known to lead to a Hopf algebra bicrossproduct
[3,7,12]. Here we also carry out the corresponding procedure for almost groups
and almost Hopf algebras.
 
This is not the only generalisation of the idea of a Hopf algebra. In 
[4] there are axioms for weak $C^{*}$-Hopf algebras, but in this case the unit and 
counit are not algebra maps, which they are in our axioms. 

Note that although we use some continuous examples of almost groups as motivation,
in the detailed proofs of the results about the algebras we shall always assume that
the almost groups are finite, or at least discrete. 

The authors would like to thank M.V.\ Lawson (Bangor) and S.\ Majid (QMW) for their assistance.

\heading{2. Almost groups}\endheading

\proclaim{Definition 2.1}
An {\bf almost group} is a set $G$ with an associative binary operation $\cdot\ ,$ a 1-1 correspondence
$i: G\longrightarrow G$ (written $g \longmapsto g^{i}$) and a set $J \subset G$ 
which is closed under the binary operation $\cdot$ and $i$. Also the following properties are satisfied:
\newline
\noindent
{\bf 2.1.1)} $(g\ h)^{i}=h^{i}\ g^{i}$, for all $g,\ h \in G$.
\newline
\noindent
{\bf 2.1.2)} For all $g \in G$, and for all $j \in J$, $j\ g=g\ j$.
\newline
\noindent
{\bf 2.1.3)} For all $g \in G$, $g\ g^{i}=g^{i}\ g \in J$.
\newline
\noindent
{\bf 2.1.4)} $(g^{i})^{i}=g,$ for all $g \in G$.
\endproclaim
 
\proclaim{Example 2.2} In the case where $J=\{e\}$, where $ge=g$
 for all $g\in G$, we just get a group.\endproclaim

\proclaim{Example 2.3}
Take $G$ to be the set of meromorphic functions from $\Bbb C_{\infty}$ to $ GL_{n}(\Bbb C)$
which are unitary on the real axis, normalised to the identity
at infinity, and have all poles in the upper half plane. 
All such loops can be factored 
as a product of functions of the form 1.1. We define the $i$ operation on 1.1 by 
$$\phi^{i}(\lambda)=\frac{\lambda - \overline\alpha}{\lambda - \alpha} P^{\bot}+ P.$$
Here $G$  is an almost group when $J$ consists of all meromorphic complex valued 
functions times the identity matrix (the `vacuum loops'),
and the binary operation is the usual matrix multiplication. \endproclaim

\proclaim{Example 2.4,  Milnor's construction of the total space of the classifying bundle of a 
group [8,9,11]}
Suppose that we have an abelian group $A$ and a topological space $EA$
consisting of step 
functions :\ $[0,1)\longrightarrow A\times \Bbb N$ which are constant on the intervals
$[x_{i},x_{i+1})$ where $0=x_{0}<x_{1}<.......<x_{N}=1$, and where the $\Bbb N$ valued numbers
are strictly increasing in successive intervals.

The multiplication on $EA$ is defined by pointwise multiplication of
the abelian group elements and adding the numbers,
and the $i$ operation is defined by taking the inverse 
of the abelian group elements without any change in the 
numbers.

Then we can say that $EA$ is an almost group if we take $J$ to consist of steps all taking value
the identity element in $A$, 
but with different numbers.\endproclaim

\proclaim{Example 2.5} Suppose that $G=\{a,b,c\},$ and $J=\{a,b\}$. We define the multiplication and the $i$
operation by $x\cdot y=a$ and $ x^{i}=x$ for all $x,y\in G$.

Alternatively we could take the case where
$a\cdot x=x\cdot a=x$ for all $x\in G$, and the other products are equal to $b$,
and the $i$ operation is as defined before.\endproclaim

\proclaim{Example 2.6} Consider $G=A\times A$ where $A$ is an abelian group,
 with multiplication $(a,b)\cdot (c,d)=(ac,bd)$, the $i$
operation $ (a,b)^{i}=(b,a)$ and
 $J=\{(a,a):a\in A\}$. 

\endproclaim

\proclaim{Example 2.7} A Clifford semigroup is an example of an
almost group, where J consists of idempotents [6].

\endproclaim

\heading{3. Almost Hopf algebras} \endheading
Now we are in the position to give a definition for an almost Hopf algebra $H$ which 
has the same rules for Hopf algebra $H$ except $\epsilon:H\longrightarrow H_J$ 
and $\eta:H_J\longrightarrow H$ where $H_J\subset H$.

\proclaim{Definition 3.1}
An almost Hopf algebra $(H\ , +\ , \cdot\ , \eta\ , \Delta\ , \epsilon\ , S\ ; H_J)$
is an associative algebra $H$ with the following additional structure:

\noindent 
{\bf 3.1.0)}\quad $H_J$ is a commutative associative algebra.

\noindent 
{\bf 3.1.1)}\quad A coassociative comultiplication map $\Delta:H\longrightarrow H\otimes H$.

\noindent 
{\bf 3.1.2)}\quad A map $\epsilon:H\longrightarrow H_J$ satisfying
$
(\epsilon\otimes id)\Delta(h)=\tau(id\otimes \epsilon)\Delta(h)
$ for all $h\in H$, where $\tau(h\otimes j)=j\otimes h$.

\noindent 
{\bf 3.1.3)}\quad A  map $\eta:H_J\longrightarrow H$ satisfying 
$
\cdot(\eta\otimes id)=\cdot(id\otimes \eta)\tau:H_J\otimes H\to H
$.

\noindent 
{\bf 3.1.4)}\quad A linear antipode map 
$S:H\longrightarrow H$ obeying
$\cdot(S\otimes id)\circ \Delta(h)=\cdot(id\otimes S)\circ \Delta(h)=\eta\circ \epsilon(h)$
 for all $h\in H$.

\noindent 
{\bf 3.1.5)}\quad The maps
 $\Delta$, $\epsilon$, and $\eta$ are algebra maps.

\endproclaim

If $G$ were a group, then its group algebra $kG$ and its function algebra $k(G)$ (where
$k$ is a field) would be Hopf algebras. We need to check that the same is true of almost
groups and almost Hopf algebras. 

\proclaim{Example 3.2} 
Take a basis of $k(G)$ consisting of elements of the form $\delta_{x}$ for $x\in G$,
and a basis for $k(G)_J=k(J)$ of the form $\delta_{j}$ for $j\in J$. 
The operations are defined as
$$\delta_{x}\cdot\delta_{y}=\delta_{x}\ \delta_{x,y}\ ,\qquad\qquad
\Delta(\delta_{x})=\sum_{x=yz} \delta_{y}\otimes\delta_{z}\ ,$$
$$S(\delta_{x})=\delta_{x^{i}}\ ,\qquad
\epsilon(\delta_{x})=\cases 0,&\text{if $x\not\in J$}\\
\delta_{x},&\text{if $x\in J$}\ ,\endcases\qquad
\eta(\delta_{j})=\sum_{z\in G:j=zz^{i}} \delta_{z}\ .$$
Here the symbol $\delta_{x,y}$ is the Kroneker delta, which is one if $x=y$ and zero otherwise.
We shall now check the rules of an almost Hopf algebra.
\endproclaim
\demo{Check} It is easy to check that $\cdot$ gives an associative multiplication
and that $\Delta$ gives a coassociative comultiplication, i.e. 
$(\Delta\otimes id)\Delta(\delta_{x})=(id\otimes \Delta)\Delta(\delta_{x})$. 
For 3.1.2 we get
$$
\align
(\epsilon\otimes id)\Delta(\delta_{x})&=(\epsilon\otimes id)(\sum_{x=yz} \delta_{y}\otimes \delta_{z})
=\sum_{x=jz:\ j\in J} \delta_{j}\otimes \delta_{z}\ , \\
\tau(id\otimes\epsilon)\Delta(\delta_{x})&=\tau(id\otimes\epsilon)(\sum_{x=zy} \delta_{z}\otimes \delta_{y})
=\tau(\sum_{x=zj:\ j\in J} \delta_{z}\otimes \delta_{j})\ ,
\endalign
$$
which are the same as $zj=jz$ for all $j\in J$. For 3.1.3,
$$
\align
\cdot(\eta\otimes id)(\delta_{j}\otimes \delta_{x})&=\cdot(\sum_{j=zz^{i}} \delta_{z}\otimes \delta_{x})
=\sum_{j=zz^{i}} \delta_{z}\ \delta_{z,x}\ ,
\\
\cdot(id\otimes \eta)\tau(\delta_{j}\otimes \delta_{x})&=\cdot(id\otimes \eta)(\delta_{x}\otimes \delta_{j})
=\cdot(\delta_{x}\otimes \sum_{j=zz^{i}} \delta_{z})
=\sum_{j=zz^{i}} \delta_{x} \delta_{x,z}\ .
\endalign
$$
For 3.1.4,
$$
\align
\cdot(S\otimes id)\Delta(\delta_{x})&=\cdot(S\otimes id)(\sum_{x=yz} \delta_{y}\otimes \delta_{z})
=\cdot(\sum_{x=yz} \delta_{y^{i}}\otimes \delta_{z})
=\sum_{x=yy^{i}} \delta_{y^{i}}\ ,\\
\cdot(id\otimes S)\Delta(\delta_{x})&=\cdot(id\otimes S)(\sum_{x=yz} \delta_{y}\otimes \delta_{z})
=\cdot(\sum_{x=yz} \delta_{y}\otimes \delta_{z^{i}})
=\sum_{x=yy^{i}} \delta_{y}\ ,
\endalign
$$
where we have used the fact that if $y=z^i$ then $y^i=(z^i)^i=z$ by 2.1.4. But these two expressions 
are the same, as can be seen by putting $w=y^i$ in the second and noting that $yy^i=y^iy$ by 2.1.3.
Now note that the expressions give zero unless $x\in J$, as $yy^i\in J$, and then they have value
$\eta\left(\epsilon(\delta_{x})\right)$. It is left to the reader to check that 
$\Delta$, $\epsilon$, and $\eta$ are algebra maps.\quad$\square$

\enddemo

\proclaim{Example 3.3} The almost group algebra has $H=kG$ and $H_J=kJ$, with multiplication given by the 
usual linear extension of the almost group multiplication. The operations are (for $x\in G$ and
$j\in J$) $\Delta(x)=x\otimes x$, $S(x)=x^{i}$, $\epsilon(x)=x\ x^{i}$ and $\eta(j)=j$.
Then the rules for an almost Hopf algebra are satisfied.\endproclaim
\demo{Check} 
 It is easy to check that $\cdot$ gives an associative multiplication
and that $\Delta$ gives a coassociative comultiplication. For 3.1.2;
$$\align
(\epsilon\otimes id)\Delta(x)&=(\epsilon\otimes id)(x\otimes x)
=x\ x^{i}\otimes x\ ,  \\
\tau(id\otimes\epsilon)\Delta(x)&=\tau(id\otimes\epsilon)(x\otimes x)
=\tau(x\otimes x\, x^{i})
=x\, x^{i}\otimes x\ .
\endalign$$
For 3.1.3;
$$\align
\cdot(\eta\otimes id)(j\otimes x)&=\cdot(j\otimes x)
=j x\ , \\
\cdot(id\otimes\eta)\tau(j\otimes x)&=\cdot(id\otimes\eta)(x\otimes j)
=\cdot(x\otimes j)
=xj=j x\ ,
\endalign$$
where we have used 2.1.2. For 3.1.4;
$$
\align
\cdot(S\otimes id)\Delta(x)&=\cdot(S\otimes id)(x\otimes x)
=\cdot(x^{i}\otimes x)
=x^{i} x\ ,  \\
\cdot(id\otimes S)\Delta(x)&=\cdot(id\otimes S)(x\otimes x)
=\cdot(x\otimes x^{i})=xx^i
=x^{i} x\ ,
\endalign
$$
where we have used 2.1.3, so both of these expressions are equal to 
$\eta\left(\epsilon(x)\right)$. It is easy to see that 
$\Delta$ and $\eta$ are algebra maps. For $\epsilon$, 
$$
\epsilon(x)\epsilon(y)=xx^iyy^i=xyy^ix^i=xy(xy)^i=\epsilon(xy)\ ,
$$
where we have used 2.1.2 and 2.1.1.\quad$\square$

\enddemo

\proclaim{Proposition 3.4} In examples 3.2 and 3.3,
$S$ is an antialgebra map, i.e.\ 
$S(h\,h')=S(h')\, S(h)$.
\endproclaim
\demo{Proof} This is immediate in both cases.\quad$\square$
\enddemo

\heading{4. Matched pairs and doublecross products} \endheading
\proclaim{Definition 4.1} Matched pairs of almost groups\endproclaim
Suppose that $(G,J_{G})$ and $(M,J_{M})$ are almost groups. 
Now take functions $\triangleright : M \times G\longrightarrow G$, and $\triangleleft : 
M \times G\longrightarrow M$,
which obey the following rules, 
for all $s,t\in M$, and $u,v\in G$: 
$$\align
s\triangleleft uv=(s\triangleleft u)\triangleleft v \quad,&\quad 
st\triangleleft u=\left(s\triangleleft (t\triangleright u)\right) (t\triangleleft u)\ , \\
st\triangleright u=s\triangleright (t\triangleright u) \quad,&\quad 
s\triangleright uv=(s\triangleright u) \left((s\triangleleft u)\triangleright v\right)\ .
\endalign$$
With the $i$ operations we require
$$(s\triangleleft u)^{i} \triangleright (s\triangleright u)^{i}=u^{i}\quad,\quad
(s\triangleleft u)^{i} \triangleleft (s\triangleright u)^{i}=s^{i}\ ,$$
and also, for all $j\in J_{G}$ or $J_{M}$, we have
$$j\triangleright u=u\ ,\quad
\quad j\triangleleft u=j\ ,\quad
s\triangleright j=j\ ,\quad s\triangleleft j=s\ .$$

\proclaim{Proposition 4.2}
Given a matched pair $(G,J_{G})$ and $(M,J_{M})$ of almost groups we can construct a 
 doublecross product almost group $G\bowtie M$, which consists of the set $G\times M$
with binary operation 
$(u,s)(v,t)=\left(u(s\triangleright v),(s\triangleleft v)t\right)$
and $(u,s)^{i}=(s^{i}\triangleright u^{i},s^{i}\triangleleft u^{i})$. Finally we set
 $J_{G\bowtie M}=J_{G}\times J_{M}$.  \endproclaim
\demo{Proof} The proof that the multiplication is associative is tedious but straight forward.
For 2.1.1,
$$\align
((u,s)(v,t))^i&=\left(u(s\triangleright v),(s\triangleleft v)t\right)^i=
(t^i(s\triangleleft v)^i\triangleright (s\triangleright v)^iu^i,
t^i(s\triangleleft v)^i\triangleleft (s\triangleright v)^iu^i) \\
&=(t^i\triangleright v^i(s^i\triangleright u^i),(t^i\triangleleft v^i)s^i\triangleleft u^i)\ , \\
(v,t)^i(u,s)^i&=(t^{i}\triangleright v^{i},t^{i}\triangleleft v^{i})
(s^{i}\triangleright u^{i},s^{i}\triangleleft u^{i})\ ,
\endalign$$
and these can be seen to be the same after the application of the
product rule. For 2.1.2, given $j\in J_G$ and $n\i J_M$,
$$\align
(j,n)(u,s)&=\left(j(n\triangleright u),(n\triangleleft u)s\right)=(ju,ns)\ ,\\
(u,s)(j,n)&=\left(u(s\triangleright j),(s\triangleleft j)n\right)=(uj,sn)\ .
\endalign$$
For 2.1.3, 
$$\align
(u,s)(u,s)^i&=(u,s)(s^{i}\triangleright u^{i},s^{i}\triangleleft u^{i})=
(uu^i,(s\triangleleft(s^{i}\triangleright u^{i}))(s^{i}\triangleleft u^{i}))\\
&=(uu^i,ss^{i}\triangleleft u^{i})=(uu^i,ss^{i})\ , \\
(u,s)^i(u,s)&=(s^{i}\triangleright u^{i},s^{i}\triangleleft u^{i})(u,s)=
((s^{i}\triangleright u^{i})((s^{i}\triangleleft u^{i})\triangleright u),s^is)\\ &=
(s^{i}\triangleright u^{i}u,s^is)=( u^{i}u,s^is)\ .
\endalign$$
For 2.1.4,
$$\align
((u,s)^i)^i&=(s^{i}\triangleright u^{i},s^{i}\triangleleft u^{i})^i=
((s^{i}\triangleleft u^{i})^{i}\triangleright (s^{i}\triangleright u^{i})^i,
(s^{i}\triangleleft u^{i})^{i}\triangleleft (s^{i}\triangleright u^{i})^i)\\
&=((u^i)^i,(s^i)^i)=(u,s)\ .\square
\endalign$$
\enddemo

\proclaim{Definition 4.3} Bicrossproduct almost Hopf algebras\endproclaim
Now from the matched pair $(G,J_{G})$ and $(M,J_{M})$
 we define an almost Hopf algebra $H=kM\subord k(G)$ with basis $s\otimes 
\delta_{u}$ where $s\in M$ and $u\in G$. We take $H_J=kJ_{M}\otimes k(J_{G})\subset H$. 
Here $kM$ is the almost group almost Hopf algebra of $M$ and 
$k(G)$ is the almost Hopf algebra of functions on $G$. Explicit formulae for $kM\subord k(G)$ are as follows:
$$
(s\otimes \delta_{u})(t\otimes \delta_{v})=\delta_{u,t\triangleright v}(st\otimes \delta_{v}),\quad\quad       
\Delta(s\otimes \delta_{u})=\sum_{xy=u} s\otimes \delta_{x}\otimes s\triangleleft x \otimes \delta_{y}
$$
$$
S(s\otimes \delta_{u})=(s\triangleleft u)^{i}\otimes \delta_{(s\triangleright u)^{i}},
 \quad\quad \epsilon(s\otimes \delta_{u})=
\cases 0,&\text{if $u\not\in J_{G}$}\\
ss^{i}\otimes \delta_{u},&\text{if $u\in J_{G}$}.\endcases
$$
$$
\eta(j\otimes \delta_{n})=
\sum_{n=zz^{i}} j\otimes \delta_{z},\qquad\text{for all $j\in J_M,\ n\in J_G$}.
$$
Now we would like to check the rules for an almost Hopf algebra, 
but first we need to prove certain results:

\proclaim{Proposition 4.4}  For all $s\in M$ and $w\in G$,
$(s\triangleleft w)^{i}(s\triangleleft w)=ss^{i}=s^{i}s$.
\endproclaim
\demo{Proof}
From the rules for a matched pair,
$$
\align
(s\triangleleft w)^{i}(s\triangleleft w)&=
(s\triangleleft w)^{i}(s\triangleleft w) \triangleleft w^{i} \\
 &=
\left((s\triangleleft w)^{i}\triangleleft \left((s\triangleleft w) 
\triangleright w^{i}\right)\right)\left((s\triangleleft w)\triangleleft w^{i}\right) \\
&=\left((s\triangleleft w)^{i}\triangleleft \left((s\triangleleft w) 
\triangleright w^{i}\right)\right)(s\triangleleft ww^{i}) \\
&=\left((s\triangleleft w)^{i}\triangleleft \left((s\triangleleft w) 
\triangleright w^{i}\right)\right) s\ .\endalign$$
Now we know that, (where $(s\triangleright w)^{i}(s\triangleright w)=j$)
$$
\align
s\triangleright ww^{i}&=
(s\triangleright w)\left((s\triangleleft w) \triangleright w^{i}\right)=ww^{i}\\
(s\triangleright w)^{i}(s\triangleright w)\left((s\triangleleft w) \triangleright w^{i}\right)&=
(s\triangleright w)^{i}ww^{i}\\
(s\triangleleft w)^{i}\triangleleft j\left((s\triangleleft w) \triangleright w^{i}\right)&=
(s\triangleleft w)^{i}\triangleleft (s\triangleright w)^{i}ww^{i}\\
(s\triangleleft w)^{i}\triangleleft \left((s\triangleleft w) \triangleright w^{i}\right)&=
\left((s\triangleleft w)^{i}\triangleleft (s\triangleright w)^{i}\right)\triangleleft ww^{i}\\
&=
s^{i}\triangleleft ww^{i}=s^{i}\ .
\square\endalign$$
\enddemo

\proclaim{Proposition 4.5}  For all $s\in M$ and $w\in G$,
 $(s\triangleright w)(s\triangleright w)^{i}=ww^{i}=w^{i}w$.
\endproclaim
\demo{Proof} From the rules for a matched pair,
$$
\align
(s\triangleright w)(s\triangleright w)^{i}&=
s^{i}\triangleright (s\triangleright w)(s\triangleright w)^{i} \\
&=
\left(s^{i}\triangleright (s\triangleright w)\right)\left(\left(s^{i}\triangleleft (s\triangleright w)\right)
\triangleright (s\triangleright w)^{i}\right)\\
&=(s^{i}s\triangleright w)\left(\left(s^{i}\triangleleft (s\triangleright w)\right)
\triangleright (s\triangleright w)^{i}\right)\\
&=w\left(\left(s^{i}\triangleleft (s\triangleright w)\right)\triangleright (s\triangleright w)^{i}\right)\ .
\endalign$$
Now we know that, (where $(s\triangleleft w)(s\triangleleft w)^{i}=j$)
$$
\align
s^{i}s\triangleleft w&=
\left(s^{i}\triangleleft (s\triangleright w)\right) (s\triangleleft w)=s^{i}s\\
\left(s^{i}\triangleleft (s\triangleright w)\right) (s\triangleleft w)(s\triangleleft w)^{i}&=
s^{i}s(s\triangleleft w)^{i}\\
\left(s^{i}\triangleleft (s\triangleright w)\right) j\triangleright (s\triangleright w)^{i} &=
s^{i}s(s\triangleleft w)^{i}\triangleright (s\triangleright w)^{i}\\
\left(s^{i}\triangleleft (s\triangleright w)\right)\triangleright (s\triangleright w)^{i} &=
s^{i}s\triangleright \left((s\triangleleft w)^{i}\triangleright (s\triangleright w)^{i}\right)\\
&=s^{i}s\triangleright w^{i}= w^{i}\ .\square
\endalign
$$
\enddemo

\demo{Check that the construction in 4.3  gives an 
 almost Hopf algebra} It is fairly standard to check that the product
is associative and that the coproduct
is coassociative. For 3.1.2, where $j\in J_G$,
$$
\align
(\epsilon\otimes id)\Delta(s\otimes \delta_{u})&=(\epsilon\otimes id)(\sum_{u=xy} s\otimes \delta_{x}
\otimes s\triangleleft x\otimes \delta_{y})\\
&=\sum_{u=jy} ss^{i}\otimes \delta_{j}
\otimes s\triangleleft j\otimes \delta_{y}
=\sum_{u=jy} ss^{i}\otimes \delta_{j}
\otimes s\otimes \delta_{y}\ ,  \\
\tau(id\otimes \epsilon)\Delta(s\otimes \delta_{u})&=\tau(id\otimes \epsilon)(\sum_{u=yx} s\otimes \delta_{y}
\otimes s\triangleleft y\otimes \delta_{x})\\
&=\tau(\sum_{u=yj} s\otimes \delta_{y}
\otimes (s\triangleleft y)(s\triangleleft y)^{i}\otimes \delta_{j})\\
&=\tau(\sum_{u=jy} s\otimes \delta_{y}
\otimes ss^{i}\otimes \delta_{j})=\sum_{u=jy} ss^{i}\otimes \delta_{j}
\otimes s\otimes \delta_{y}\ .
\endalign
$$
For 3.1.3, where $j\in J_M$ and $n\in J_G$,
$$
\align
\cdot(\eta\otimes id)&\left((j\otimes \delta_{n})\otimes (s\otimes \delta_{u})\right)=
\cdot(\sum_{n=zz^{i}} (j\otimes \delta_{z})\otimes (s\otimes \delta_{u}))\\
&=\sum_{n=zz^{i}} \delta_{z, s\triangleright u}\ js\otimes \delta_{u}
=\sum_{n=uu^{i}} js\otimes \delta_{u}\ ,\\
\cdot(id\otimes \eta)&\tau\left((j\otimes \delta_{n})\otimes (s\otimes \delta_{u})\right)=
\cdot(id\otimes \eta)((s\otimes \delta_{u})\otimes (j\otimes \delta_{n}) )\\
&=\cdot((s\otimes \delta_{u})\otimes (\sum_{n=zz^{i}} j\otimes \delta_{z}))
=\sum_{n=zz^{i}} \delta_{u, j\triangleright z} \ sj\otimes \delta_{z} \\
&=\sum_{n=zz^{i}} \delta_{u,z} \ js\otimes \delta_{z}
=\sum_{n=uu^{i}} js\otimes \delta_{u}
\endalign
$$
For 3.1.4, 
$$ 
\align
\cdot(S\otimes id)\Delta(s\otimes \delta_{u})&=
\cdot(S\otimes id)(\sum_{xy=u} s\otimes \delta_{x}\otimes (s\triangleleft x)\otimes \delta_{y}) \\
&=\sum_{xy=u}\left((s\triangleleft x)^{i}\otimes \delta_{(s\triangleright x)^{i}}\right)
\left((s\triangleleft x)\otimes \delta_{y}\right) \\
&=\sum_{xy=u} \delta_{(s\triangleright x)^{i},(s\triangleleft x)\triangleright y} 
\ (s\triangleleft x)^{i}
(s\triangleleft x)\otimes \delta_{y}\tag*\\
&=\sum_{xy=u} \delta_{(s\triangleright x)^{i},(s\triangleleft x)\triangleright y} 
\  s^{i}s\otimes \delta_{y}\tag**\\
&=\sum_{xy=u} \delta_{x^{i}, y} 
\  s^{i}s\otimes \delta_{y}=\sum_{y^iy=u} 
  s^{i}s\otimes \delta_{y}\ ,
\endalign$$
where we have used 4.4 on (*), and applied $(s\triangleleft x)^{i}\triangleright$ to both elements in
$\delta_{(s\triangleright x)^{i},(s\triangleleft x)\triangleright y} $ in (**). 
$$\align
\cdot(id\otimes S)\Delta(s\otimes \delta_{u})&=
\cdot(id\otimes S)(\sum_{xy=u} s\otimes \delta_{x}\otimes (s\triangleleft x)\otimes \delta_{y})\\
&=\sum_{xy=u}
(s\otimes \delta_{x})\left(\left(s\triangleleft x y\right)^{i}\otimes 
\delta_{\left((s\triangleleft x)\triangleright y\right)^{i}}\right) \\
&=\sum_{xy=u} 
\delta_{x,\left((s\triangleleft x)\triangleleft y\right)^{i}\triangleright 
\left((s\triangleleft x)\triangleright y\right)^{i}}
\, s\left(s\triangleleft x y\right)^{i}\otimes 
\delta_{\left((s\triangleleft x)\triangleright y\right)^{i}} \\
&=\sum_{xy=u} \delta_{x,y^{i}}\,s(s\triangleleft xy)^{i}\otimes 
\delta_{\left((s\triangleleft x)\triangleright y\right)^{i}}\tag ***
\\
&=\sum_{xy=u} \delta_{x,y^{i}}\,ss^{i}\otimes 
\delta_{\left((s\triangleleft x)\triangleright y\right)^{i}}=
\sum_{y^iy=u}\,ss^{i}\otimes 
\delta_{\left((s\triangleleft y^i)\triangleright y\right)^{i}}\ .
\endalign
$$
In (***)
we have used the fact that $u=xy=y^iy\in J_G$. Now change variable in the sum from
$y$ to $z^i=(s\triangleleft y^i)\triangleright y$.   Then
$((s\triangleleft y^i)\triangleleft y)^i\triangleright((s\triangleleft y^i)\triangleright y)^i=y^i$, 
so $s^i\triangleright z=y^i$. Then the condition on the summation is
$u=yy^i=(s^i\triangleright z)^i(s^i\triangleright z)=z^iz$
This shows that the
sums are the same, and that they have value
$\eta\left(\epsilon(s\otimes \delta_{u})\right)$.

To show that $\Delta$ is an algebra map,
$$ 
\align
\Delta \left((s\otimes \delta_{u})(t\otimes \delta_{v})\right)&=
\Delta\left(\delta_{u,t\triangleright v}\,st\otimes \delta_{v}\right) 
=\sum_{xy=v} \delta_{u,t\triangleright v}\,
st\otimes \delta_{x}\otimes st\triangleleft x \otimes \delta_{y} \\
&=\sum_{xy=v=t^{i}\triangleright u} st\otimes \delta_{x}\otimes st\triangleleft x \otimes \delta_{y}\ ,\\
\Delta (s\otimes \delta_{u})\Delta  (t\otimes \delta_{v})&=
(\sum_{xy=u} s\otimes \delta_{x}\otimes s\triangleleft x \otimes \delta_{y})
(\sum_{x_{1}y_{1}=v} t\otimes \delta_{x_{1}}\otimes 
t\triangleleft x_{1} \otimes \delta_{y_{1}}) \\
&=\sum_{xy=u,x_{1}y_{1}=v} (s\otimes \delta_{x})(t\otimes 
\delta_{x_{1}})\otimes (s\triangleleft x \otimes \delta_{y})
(t\triangleleft x_{1} \otimes \delta_{y_{1}}) \\
&=\sum_{xy=u,x_{1}y_{1}=v} \delta_{x,t\triangleright x_{1}}
\delta_{y,(t\triangleleft x_{1})\triangleright y_{1}}\,st\otimes \delta_{x_{1}}\otimes
(s\triangleleft x)(t\triangleleft x_{1})\otimes \delta_{y_{1}} \\
&=\sum_{xy=u,x_{1}y_{1}=v} \delta_{x,t\triangleright x_{1}}
\delta_{y,(t\triangleleft x_{1})\triangleright y_{1}}
\,st\otimes \delta_{x_{1}}\otimes st\triangleleft x_{1}\otimes \delta_{y_{1}} \\
&=\sum_{x_{1}y_{1}=v=t^{i}\triangleright u} st\otimes \delta_{x_{1}}\otimes 
st\triangleleft x_{1}\otimes \delta_{y_{1}}\ .
\endalign$$
To show that $\epsilon$ is an algebra map, if $u,v\in J_{G}$,
$$ 
\align
\epsilon (s\otimes \delta_{u})\epsilon (t\otimes \delta_{v})&=(ss^{i}\otimes \delta_{u})
(tt^{i}\otimes \delta_{v})
=\delta_{u,tt^{i}\triangleright v} \,ss^{i}tt^{i}\otimes \delta_{v}  \\
&=\delta_{u,v} \,stt^{i}s^{i}\otimes \delta_{v}
=\delta_{u,v} (st)(st)^{i}\otimes \delta_{v}\ , \\
\epsilon\left((s\otimes \delta_{u})(t\otimes \delta_{v})\right)&=
\epsilon \left(\delta_{u,t\triangleright v} \,st\otimes \delta_{v}\right) 
=\delta_{u,v} (st)(st)^{i}\otimes \delta_{v}\ .\endalign
$$
If  $u\notin J_{G}$ or $v\notin J_{G}$ then both expressions give zero. 
To show that $\eta$ is an algebra map,
$$
\align
\eta (j\otimes \delta_{n})&\eta (j_{1}\otimes \delta_{n_{1}})=
(\sum_{zz^{i}=n} j\otimes \delta_{z})
(\sum_{xx^{i}=n_{1}} j_{1}\otimes \delta_{x}) \\
&=\sum_{zz^{i}=n,xx^{i}=n_{1}} \delta_{z,j_{1}\triangleright x} \,jj_{1}\otimes \delta_{x} \\
&=\sum_{zz^{i}=n,xx^{i}=n_{1}} \delta_{z,x} \,jj_{1}\otimes \delta_{x}
=\sum_{zz^{i}=n=n_{1}} jj_{1}\otimes \delta_{z}\ ,
\\
\eta((j\otimes \delta_{n})&(j_{1}\otimes \delta_{n_{1}}))=
\eta(\delta_{n,j_{1}\triangleright n_{1}} \,jj_{1}\otimes \delta_{n_{1}})\\ 
&
=\sum_{zz^{i}=n_{1}}\delta_{n,n_{1}}\,jj_{1}\otimes \delta_{z}\ .
\endalign$$

\enddemo

\proclaim{Proposition 4.6}
$S$ reverses the order of the product and coproduct, and preserves $S$.
\endproclaim
\demo{Proof} For the product:
$$
\align
S\left((s\otimes \delta_{u})(t\otimes \delta_{v})\right)&=
S\left(\delta_{u\ ,t\triangleright v} \ st\otimes \delta_{v}\right)
=\delta_{u\ ,t\triangleright v} \ (st\triangleleft v)^{i}\otimes 
\delta_{(st\triangleright v)^{i}} \\
&=\delta_{u\ ,t\triangleright v} \ \left(\left(s\triangleleft (t\triangleright v)
\right)(t\triangleleft v)\right)^{i}\otimes 
\delta_{\left(s\triangleright (t\triangleright v)\right)^{i}} \\
&=\delta_{u\ ,t\triangleright v} \ \left((s\triangleleft u)(t\triangleleft v)\right)^{i}\otimes 
\delta_{(s\triangleright u)^{i}}\\
&=\delta_{u\ ,t\triangleright v} \ (t\triangleleft v)^{i}(s\triangleleft u)^{i}\otimes 
\delta_{(s\triangleright u)^{i}}\ ,\\
S(t\otimes \delta_{v}) S(s\otimes \delta_{u})&=\left((t\triangleleft v)^{i}\otimes  
\delta_{(t\triangleright v)^{i}}\right)
\left((s\triangleleft u)^{i}\otimes  \delta_{(s\triangleright u)^{i}}\right)\\
&=\delta_{(t\triangleright v)^{i} ,(s\triangleleft u)^{i}\triangleright (s\triangleright u)^{i}} 
\ (t\triangleleft v)^{i}(s\triangleleft u)^{i}\otimes \delta_{(s\triangleright u)^{i}}\\
&=\delta_{(t\triangleright v)^{i} ,u^{i}} 
\ (t\triangleleft v)^{i}(s\triangleleft u)^{i}\otimes \delta_{(s\triangleright u)^{i}} 
\ .
\endalign
$$
For the coproduct:
$$\align
\tau\Delta 
S(s\otimes \delta_{u})&=\tau\Delta((s\triangleleft u)^{i}\otimes \delta_{(s\triangleright u)^{i}})
\\ &=\sum_{xy=(s\triangleright u)^{i}}
(s\triangleleft u)^{i}\triangleleft x\otimes \delta_{y}
\otimes
(s\triangleleft u)^{i}\otimes \delta_{x} \ ,\\
(S\otimes S)\Delta (s\otimes \delta_{u})&=\sum_{vw=u}
(s\triangleleft v)^{i}\otimes \delta_{(s\triangleright v)^{i}}
\otimes
(s\triangleleft u)^{i}\otimes \delta_{((s\triangleleft v)\triangleright w)^{i}}\ .
\endalign$$
Now we set $x=((s\triangleleft v)\triangleright w)^{i}$ and $y=(s\triangleright v)^{i}$
in the first expression, and observe that we get the second expression.
As for the antipode, $S$ commutes with itself. \quad$\square$
\enddemo

\proclaim{Proposition 4.7}
We have $\epsilon S=S_J\epsilon$ and $S\eta=\eta S_J$, where $S_J$ is the restriction of $S$
to $H_J$, which is just $S_J(j\otimes \delta_n)=j^i\otimes \delta_{n^i}$.
\endproclaim
\demo{Proof} This is fairly simple, and is left to the reader.\quad$\square$
\enddemo

\heading{5. The meromorphic loop group}\endheading
In this section we continue with the meromorphic loops introduced in example 2.3.
Any invertible
 meromorphic function $\phi:\Bbb C_{\infty}\longrightarrow M_{n}$ which is unitary on the real axis can be
written as a constant matrix times a product of factors of the form 
$$
\Phi_{\alpha,P}(\lambda)= P^{\bot} + \frac{\lambda - \overline \alpha}{\lambda - 
\alpha} P\ ,$$
where $\alpha\in \Bbb C\setminus \Bbb R$ and $P$ is a self-adjoint projection in $M_{n}$. 
Define the sets
$$\align
G&=\{\text{$\phi$ : $\phi(\infty)=1$ and $\phi(\lambda)$ has no singularities for im$(\lambda)> 0$} \}\ ,\\
J_{G}&=\{\text{$\phi$ : $\phi(\infty)=1$ and $\phi(\lambda)
$ has no singularities for im$(\lambda)> 0$, and}\\
&\qquad \text{$\phi$ is a scalar function times the identity matrix}\} \ ,\\
M&=\{\text{$\phi$ : $\phi(\infty)=1$ and $\phi(\lambda)
$ has no singularities for im$(\lambda)< 0$} \}\ ,\\
J_{M}&=\{\text{$\phi$ : $\phi(\infty)=1$ and $\phi(\lambda)$ has no singularities for im$(\lambda)< 0$, and}\\
&\qquad \text{$\phi$ is a scalar function times the identity matrix}\} \ .
\endalign$$
The normalisation $\phi(\infty)=1$ just means that we can forget about the constant factor.
Define the $i$ operation by 
$$
\Phi_{\alpha,P}^i(\lambda)=\frac{\lambda -  \overline \alpha}{\lambda - 
\alpha}  P^{\bot}+P\ ,
$$ 
and extend this to products of basic loops by reversing order, i.e.\ $(\Phi\Psi)^i=\Psi^i\Phi^i$. 
It is not too difficult to show that $(G,J_{G}),$ and $(M,J_{M})$ are almost groups,
with the usual matrix multiplication.

\proclaim{Definition 5.1}
We define the actions 
$\triangleright$ and $\triangleleft$ by reversal of order
of multiplication, i.e.\ for $s\in M$ and $u\in G$
choose $s\triangleright u\in G$ and $s\triangleleft u\in M$ so that
$su=(s\triangleright u)(s\triangleleft u)$. Here we must issue a warning; there is no uniqueness
 of factorisation.  
To factor a meromorphic loop $\phi$ we can use the procedure in [2] to write $\phi$ as a product of 
basic loops, choosing the lower half plane poles first, $\phi=su$. 
There are other possible factorisations of the form $\phi=s'u'$, where $s'=st\in M$
and $u'=t^{-1}u\in G$, all we have to do is to take 
$t$ with all poles in the upper half plane, so $t^{-1}$ has all poles in the lower half plane. 
This occurs because $G$ and $M$ are not groups, as they are not closed
under the inverse operation.
To get round this, we always choose a factorisation with the minimum number of
basic factors. It is also possible to have ambiguities in the factorisation where
poles coincide or are at complex conjugate positions (as noted in the proof
of the following proposition). Strictly we should restrict our results
on actions to the
dense open set of loops which have no multiple poles or poles at
complex conjugate positions. We shall assume this for the rest of the section (with
the exception of the next proposition).
\endproclaim

We can calculate the actions on the basic factors by the next result. The actions
on products of basic factors are calculated by successive reversals of factors,
a procedure which does not increase the number of factors. In fact, $s\triangleright u$
has exactly the same pole positions as $u$, and $s\triangleleft u$
has exactly the same pole positions as $s$.

\proclaim{Proposition 5.2}
Suppose
$\theta_{\alpha}(\lambda)=(\lambda - \overline\alpha)/(\lambda - \alpha)$ and 
$\theta_{\beta}(\lambda)=(\lambda - \overline\beta)/(\lambda - \beta)$,
where $\alpha$ and $\beta$ are in different half planes (in 
particular $\alpha\neq\beta$).  Then
$$(P_{1}^{\bot}+\theta_{\alpha} P_{1})(P_{2}^{\bot}+\theta_{\beta} P_{2})=
(P_{3}^{\bot}+\theta_{\beta} P_{3})(P_{4}^{\bot}+\theta_{\alpha} P_{4})$$
where (if we put $V_i$ to be the image of the projection $P_i$)
$$
V_{3}=(P_{1}^{\bot}+\theta_{\alpha} (\beta) P_{1})V_{2}\quad \text{ and }\quad
V_{4}=(P_{3}^{\bot}+\theta_{\beta}^{-1} (\alpha) P_{3})V_{1}\ ,
$$
if $\beta\neq\overline\alpha$, and if $\beta=\overline\alpha$ we get
$
P_{3}=1-P_{1}$ and $
P_{4}=1-P_{2}$.
\endproclaim
\demo{Proof}
We know $P_{1}$ and $P_{2}$, and we want to get  $P_{3}$ and $P_{4}$. If $\alpha\neq
\overline\beta$, we have
$
\left(P_{1}^{\bot}+\theta_{\alpha} (\beta) P_{1}\right)V_{2}=V_{3}$ and $
\left(P_{3}^{\bot}+\theta_{\beta} (\alpha) P_{3}\right)V_{4}=V_{1}
$,
which implies
$V_{4}=(P_{3}^{\bot}+\theta_{\beta}^{-1} (\alpha) P_{3})V_{1}$.
But if $\beta=\overline\alpha$ there is a problem,
 because $\theta_{\beta}(\alpha)$ is not invertible. 
If $\beta=\overline\alpha$, we know that
$\theta_{\alpha}(\lambda)=1/\theta_{\beta}(\lambda)$. Then setting
$z=\theta_{\alpha}(\lambda)$,
 we can write the factorisation as
$$
(P_{1}^{\bot}+zP_{1})(P_{2}^{\bot}+\frac{1}{z} P_{2})=
(P_{3}^{\bot}+\frac{1}{z} P_{3})(P_{4}^{\bot}+z P_{4})\ ,
$$
which can be rearranged to give
$$
(P_{2}^{\bot}+\frac{1}{z} P_{2})(P_{4}^{\bot}+\frac{1}{z}P_{4})=
(P_{1}^{\bot}+\frac{1}{z}P_{1})(P_{3}^{\bot}+\frac{1}{z} P_{3})\ .
$$
By separating powers of $z$ we get $P_4=P_3+P_1-P_2$ and $(P_1-P_2)P_3=P_2(P_1-P_2)$. 
In the case where $P_1-P_2$ is invertible, we can define $P_3$ as the unique solution
to $(P_1-P_2)P_3=P_2(P_1-P_2)$, and this will then give a unique value of $P_4$. From substituting
in the equation we see that these unique solutions are $P_{3}=1-P_{1}$ and $
P_{4}=1-P_{2}$. To preserve continuity, we
 will define these to be the actions even if $P_1-P_2$ is not invertible.
\quad$\square$
\enddemo

\proclaim{Proposition 5.3} The meromorphic loop
almost groups $(G,J_G)$ and $(M,J_M)$, with the actions
and $i$ operation specified, form a matched pair.
\endproclaim
\demo{Proof}
Consider the associativity of the multiplication $stu$ where $s,t\in M$ and 
$u\in G$. Then,
$$
s(tu)=(st)u
=(st\triangleright u)(st\triangleleft u)
=s(t\triangleright u)(t\triangleleft u)
=\left(s\triangleright (t\triangleright u)\right)\left(s\triangleleft
 (t\triangleright u)\right)(t\triangleleft u)
.
$$
By the uniqueness of the factorisation (on the open dense subset
referred to earlier), we see that
$
s\triangleright (t\triangleright u)=st\triangleright u$ and $
\left(s\triangleleft (t\triangleright u)\right)(t\triangleleft u)=st\triangleleft u
$.
Similarly, for all $s\in M$ and $u,v\in G$, we have 
$$
(su)v=s(uv)
=(s\triangleright uv)(s\triangleleft uv)
=(s\triangleright u)(s\triangleleft u)v
=(s\triangleright u)\left((s\triangleleft u)\triangleright v\right)\left((s\triangleleft u)\triangleleft v\right)
$$
which gives 
$
(s\triangleright u)\left((s\triangleleft u)\triangleright v\right)=s\triangleright uv$ and $
(s\triangleleft u)\triangleleft v=s\triangleleft uv
$.
Also, for all $j\in J_{M}$ and $u\in G$, we have 
$
ju=uj
=(j\triangleright u)(j\triangleleft u)
$,
which gives
$
j\triangleright u=u$ and $j\triangleleft u=j
$.
Similarly, for all $j\in J_{G}$ and $s\in M$, we have 
$
sj=js
=(s\triangleright j)(s\triangleleft j)
$,
which gives
$
s\triangleright j=j$ and $s\triangleleft j=s
$.
Finally, for all $s\in M$ and $u\in G$, we have 
$$
\align
(su)^{i}=u^{i}s^{i}
=\left((s\triangleright u)(s\triangleleft u)\right)^{i}
&=(s\triangleleft u)^{i}(s\triangleright u)^{i}\\
&=\left((s\triangleleft u)^{i}\triangleright (s\triangleright u)^{i}\right)
\left((s\triangleleft u)^{i}\triangleleft (s\triangleright u)^{i}\right).\endalign
$$ 
By the uniqueness of the factorisation, we see that  
$$
(s\triangleleft u)^{i}\triangleright (s\triangleright u)^{i}=u^{i},\qquad\text{and}\qquad
(s\triangleleft u)^{i}\triangleleft (s\triangleright u)^{i}=s^{i}.\quad\square
$$

\enddemo

\heading{6. Duality}\endheading
 We take $(G,J_{G})$ and $(M,J_{M})$ to be a matched pair of almost groups.
There is a dual almost Hopf algebra $H^{'}=k(M)\subordd kG$
to $H=kM\subord k(G)$ with basis $\delta_{s}\otimes u$ where $s\in M$ 
and $u\in G$, with $J_{H'}=k(J_{M})\subordd kJ_{G}$.
The explicit formulae for this almost Hopf algebra are as follows:
$$
(\delta_{s}\otimes u)(\delta_{t}\otimes v)=\delta_{s\triangleleft u,t}(\delta_{s}\otimes uv),\quad\quad       
\Delta(\delta_{s}\otimes u)=
\sum_{a,b\in M:ab=s} \delta_{a}\otimes b\triangleright u \otimes \delta_{b} \otimes u\ ,
$$
$$
S(\delta_{s}\otimes u)=\delta_{(s\triangleleft u)^{i}} \otimes(s\triangleright u)^{i}, \quad\quad 
\epsilon(\delta_{s}\otimes u)=
\cases 0,&\text{if $s\not\in J_{M}$}\\
\delta_{s}\otimes uu^{i},&\text{if $s\in J_{M}$}\endcases\ ,
$$
$$
\eta(\delta_{j}\otimes n)=\sum_{a\in M:j=aa^{i}} \delta_{a}\otimes n
\ .
$$
The dual pairing between $H'$ and $H$ is given by
$$
\big<\delta_{s}\otimes u,t\otimes \delta_{v}\big>\,=\,\delta_{s,t}\,\delta_{u,v}\ .
$$

\proclaim{Proposition 6.1} The almost Hopf algebras $H=kM\subord k(G)$ and 
$H^{'}=k(M)\subordd kG$ are dual to each other.
\endproclaim
\demo{Proof} First we check that the counits and the units are dual to each other:
$$\align
\big<\epsilon(\delta_{s}\otimes u),j\otimes\delta_n\big>&\,=\, \delta_{s,j}\,\delta_{uu^i,n} \\
\big<\delta_{s}\otimes u,\eta(j\otimes\delta_n)\big>&\,=\, 
\big<\delta_{s}\otimes u,\sum_{n=zz^{i}} j\otimes \delta_{z}\big>\,=\,
\delta_{s,j}\,\delta_{uu^i,n}\ ,\\
\big<\delta_{j}\otimes n,\epsilon(s\otimes\delta_u)\big>&\,=\, \delta_{j,ss^i}\,\delta_{n,u} \\
\big<\eta(\delta_{j}\otimes n),s\otimes\delta_u\big>&\,=\,
\big<\sum_{j=zz^{i}} \delta_{z}\otimes n,s\otimes\delta_u\big>\,=\,\delta_{j,ss^i}\,\delta_{n,u}\ .
\endalign$$
Now we check the antipodes:
$$\align
\big<S(\delta_{s}\otimes u),t\otimes\delta_v\big>&\,=\, 
\big<\delta_{(s\triangleleft u)^{i}} \otimes(s\triangleright u)^{i},t\otimes\delta_v\big>
\,=\,\delta_{(s\triangleleft u)^{i},t}\,\delta_{(s\triangleright u)^{i},v}\ ,\\
\big<\delta_{s}\otimes u,S(t\otimes\delta_v)\big>&\,=\, \big<\delta_{s}\otimes u,
(t\triangleleft v)^{i}\otimes \delta_{(t\triangleright v)^{i}}\big>\,=\,
\delta_{s,(t\triangleleft v)^{i}}\,\delta_{u,(t\triangleright v)^{i}}\ ,
\endalign$$
and these are the same by the original definition of the actions. It is left to the reader to check
the product and coproduct, i.e.
$$\align
\big<(\delta_{s}\otimes u)(\delta_{t}\otimes v),r\otimes\delta_w\big>&\,=\, 
\big<(\delta_{s}\otimes u)\otimes(\delta_{t}\otimes v),\Delta(r\otimes\delta_w)\big>\ ,\\
\big<\delta_{s}\otimes u,(t\otimes\delta_v)(r\otimes\delta_w)\big>&\,=\, 
\big<\Delta(\delta_{s}\otimes u),(t\otimes\delta_v)\otimes(r\otimes\delta_w)\big>\ .\square
\endalign$$
\enddemo

\heading{7. The $*$ Operation}\endheading
Let us define a $ * $ operation on $H$ by
$
(s\otimes \delta_{u})^{*} =s^{i}\otimes \delta_{s\triangleright u}
$
on the basis elements, extended to a conjugate-linear map from $H$ to $H$.

\proclaim{Proposition 7.1}
The $ * $ operation reverses the order of multiplication.
\endproclaim
\demo{Proof}
$$\align
\left((s\otimes \delta_{u})(t\otimes \delta_{v})\right)^{*} &=
\left(\delta_{u,t\triangleright v}\,st\otimes \delta_{v}\right)^{*}  
=\delta_{u,t\triangleright v}\,(st)^{i}\otimes \delta_{st\triangleright v}\ ,
\\
(t\otimes \delta_{v})^{*} (s\otimes \delta_{u})^{*}& =
(t^{i}\otimes \delta_{t\triangleright v})(s^{i}\otimes \delta_{s\triangleright u}) 
=\delta_{t\triangleright v,s^{i}\triangleright (s\triangleright u)}
\,t^{i}s^{i}\otimes \delta_{s\triangleright u} \\
&=\delta_{t\triangleright v,s^{i}s\triangleright u}\,(st)^{i}\otimes \delta_{s\triangleright u} 
=\delta_{t\triangleright v,u}\,(st)^{i}\otimes \delta_{s\triangleright (t\triangleright v)}\ .
\quad\square
\endalign$$
\enddemo

\proclaim{Proposition 7.2}
The $*$ operation preserves the comultiplication.
\endproclaim
\demo{Proof}
$$\align
\Delta \left((s\otimes \delta_{u})^{*}\right)&=\Delta(s^{i}\otimes \delta_{s\triangleright u}) 
=\sum_{xy=s\triangleright u} s^{i}\otimes \delta_{x}\otimes s^{i}\triangleleft x\otimes \delta_{y}\ ,
\\
\left(\Delta (s\otimes \delta_{u})\right)^{*}&=
(\sum_{x_{1}y_{1}=u} s\otimes \delta_{x_{1}}\otimes s\triangleleft 
x_{1}\otimes \delta_{y_{1}})^{*} \\
&=\sum_{x_{1}y_{1}=u} s^{i}\otimes \delta_{s\triangleright x_{1}}\otimes (s\triangleleft x_{1})^{i}\otimes 
\delta_{(s\triangleleft x_{1})\triangleright y_{1}}.
\endalign$$
Since 
$s\triangleright u=s\triangleright x_{1}y_{1}=(s\triangleright x_{1})((s\triangleleft x_{1})
\triangleright y_{1})$, if we consider
$x=s\triangleright x_{1}$ (i.e.\ $x_{1}=s^{i}\triangleright x$) and $
y=(s\triangleleft x_{1})\triangleright y_{1}$ (i.e.\ $y_{1}=(s\triangleleft x_{1})^{i}
\triangleright y$) we see that the two sums are the same.$\square$\enddemo

\proclaim{Proposition 7.3}
The $*$ operation preserves the unit, counit and antipode.
\endproclaim
\demo{Proof} For the unit:
$$
\align
\left(\eta (j\otimes \delta_{n})\right)^{*}&=(\sum_{zz^{i}=n} j\otimes \delta_{z})^{*}
=\sum_{zz^{i}=n} j^{i}\otimes \delta_{j\triangleright z}
=\sum_{zz^{i}=n} j^{i}\otimes \delta_{z}\ ,\\
\eta\left((j\otimes \delta_{n})^{*}\right)&=\eta\left(j^{i}\otimes 
\delta_{j\triangleright n}\right)
=\sum_{zz^{i}=j\triangleright n} j^{i}\otimes \delta_{z}
=\sum_{zz^{i}=n} j^{i}\otimes \delta_{z}\ .\endalign
$$
For the counit:
$$
\align
\epsilon\left((s\otimes \delta_{u}\right)^{*})=\epsilon(s^{i}\otimes \delta_{s\triangleright u})
&=\cases 0,&\text{if $s\triangleright u\not\in J_{G}$}\\
s^{i}s\otimes \delta_{s\triangleright u},&\text{if $s\triangleright u\in J_{G}$}.\endcases\\
&=\cases 0,&\text{if $u\not\in J_{G}$}\\
s^{i}s\otimes \delta_{u},&\text{if $u\in J_{G}$}.\endcases\\
&=(\epsilon(s\otimes \delta_{u}))^*\ .\endalign
$$
For the antipode:
$$
\align
(S (s\otimes \delta_{u}))^{*}&=((s\triangleleft u)^{i}\otimes 
\delta_{(s\triangleright u)^{i}})^{*}
=((s\triangleleft u)^{i})^{i}\otimes 
\delta_{(s\triangleleft u)^{i}\triangleright (s\triangleright u)^{i}}\\
&=(s\triangleleft u)\otimes \delta_{u^{i}}\ ,\\
S((s\otimes \delta_{u})^{*})&=S(s^{i}\otimes \delta_{s\triangleright u})
=(s^{i}\triangleleft (s\triangleright u))^{i}\otimes 
\delta_{(s^{i}\triangleright (s\triangleright u))^{i}}\\
&=(s^{i}\triangleleft (s\triangleright u))^{i}\otimes 
\delta_{(s^{i}s\triangleright u)^{i}}
=(s\triangleleft u)\otimes \delta_{u^{i}}\ ,\endalign
$$ 
as $s^{i}\triangleleft (s\triangleright u)=(s\triangleleft u)^i$ (just apply 
$\triangleleft(s\triangleright u)^i$ to both sides).$\square$
 \enddemo

\heading{8. Mutually inverse matched pairs}\endheading
Here we discuss a property motivated by the meromorphic loop example discussed earlier.

\proclaim{Definition 8.1} The matched pair $(G,J_G)$ and $(M,J_M)$ is said to be
mutually inverse if  the following conditions hold:
\newline
\noindent
{\bf 8.1.1)} The doublecross product $G\bowtie M$ is a group, with inverse operation
$x\mapsto x^{-1}$. 
\newline
\noindent
{\bf 8.1.2)} For all $s\in M$, $s^{-1}\in G$ and also for all 
$u\in G$, $u^{-1}\in M$.
\newline
\noindent
{\bf 8.1.3)} The map inverse : $J_{G}\longrightarrow J_{M}$ is a 1-1 correspondence.
\newline
\noindent
{\bf 8.1.4)} For all $x\in G\bowtie M$, $(x^{-1})^{i}=(x^{i})^{-1}$.
\newline
\noindent
{\bf 8.1.5)} For all $s\in M$
and $u\in G$, $u^{-1}\triangleright s^{-1}=(s\triangleleft u)^{-1}$ and 
$u^{-1}\triangleleft s^{-1}=(s\triangleright u)^{-1}$.
\endproclaim

\proclaim{Example 8.2} 
The meromorphic loop almost groups defined in 5.1 form a
mutually inverse matched pair (with the usual caveat about densely defined actions).\endproclaim
\demo{Check} The doublecross product just consists of meromorphic loops which are unitary
on the real axis, with the usual pointwise multiplication. 
On the single pole factors the inverse is
$$
\align
\Big( P^{\bot} + \frac{\lambda - \overline \alpha}{\lambda - 
\alpha} P\Big)^{-1}&=
 P^{\bot} + \frac{\lambda - \alpha}{\lambda - \overline \alpha}
 P\ ,
\endalign
$$
so that a factor with a pole in the upper half plane has an inverse with a pole in the lower half plane,
and vice versa. It is fairly easy to check that $(x^{-1})^{i}=(x^{i})^{-1}$ from this formula.

 If we take a factorisation $su = (s\triangleright u)
(s\triangleleft u)$ (where $s\in M$ and $u\in G$), and take the inverses
of both sides we get 
$u^{-1}s^{-1} =(s\triangleleft u)^{-1}(s\triangleright u)^{-1} $. But 
$u^{-1}\in M$ and $s^{-1}\in G$, so 
$u^{-1}s^{-1} 
=(u^{-1}\triangleright s^{-1})(u^{-1}\triangleleft s^{-1})$. As both 
$u^{-1}\triangleright s^{-1}$ and $(s\triangleleft u)^{-1}$ have the same pole positions
 we see that
$u^{-1}\triangleright s^{-1}=(s\triangleleft u)^{-1}$, and similarly
$u^{-1}\triangleleft s^{-1}= (s\triangleright u)^{-1}$.\quad $\square$
\enddemo

\proclaim{Definition 8.3} In the case where we have a mutually inverse matched pair, we 
define the map $T: H=kM\subord k(G)\longrightarrow H^{'}=k(M)\subordd kG$ by
$T(s\otimes \delta_{u})=\delta_{u^{-1}}\otimes s^{-1}$, and 
$T_J:J_H\to J_{H'}$ by $T_J(j\otimes \delta_{n})=\delta_{n^{-1}}\otimes j^{-1}$.
\endproclaim

\proclaim{Proposition 8.4}
The map $T$ reverses the order of both multiplication
and comultiplication.
\endproclaim
\demo{Proof} For multiplication:
$$
\align
T((s\otimes \delta_{u})(t\otimes \delta_{v}))&=
T(\delta_{u,t\triangleright v}\ st\otimes \delta_{v}) 
=\delta_{u,t\triangleright v}\ \delta_{v^{-1}}\otimes (st)^{-1}\ , \\
T(t\otimes \delta_{v})T(s\otimes \delta_{u})&=
(\delta_{v^{-1}}\otimes t^{-1})(\delta_{u^{-1}}\otimes s^{-1}) 
=\delta_{v^{-1}\triangleleft t^{-1},u^{-1}}\ \delta_{v^{-1}}\otimes t^{-1}s^{-1} \\
&=\delta_{(t\triangleright v)^{-1},u^{-1}}\ \delta_{v^{-1}}\otimes (st)^{-1}\ .\endalign
$$
For comultiplication:
$$
\align
(T\otimes T)\left(\Delta (s\otimes \delta_{u})\right)&=
(T\otimes T)(\sum_{xy=u} s\otimes \delta_{x}\otimes s\triangleleft x\otimes \delta_{y}) \\
&=\sum_{xy=u} \delta_{x^{-1}}\otimes s^{-1}\otimes 
\delta_{y^{-1}}\otimes (s\triangleleft x)^{-1}\ ,
\\
\tau\Delta \left(T(s\otimes \delta_{u})\right)&=
\tau\Delta (\delta_{u^{-1}}\otimes s^{-1}) 
=\tau(\sum_{ab=u^{-1}} \delta_{a}\otimes b\triangleright s^{-1}\otimes 
\delta_{b}\otimes s^{-1}) \\
&=\sum_{ab=u^{-1}} 
\delta_{b}\otimes s^{-1}\otimes 
\delta_{a}\otimes b\triangleright s^{-1}\ ,
\endalign$$
which can be seen to be the same on substituting $a=y^{-1}$ and $b=x^{-1}$.\quad
$\square$
\enddemo

\proclaim{Proposition 8.5}
The map $T$ preserves the antipode  and $*$-operation, where $*$ on $H'$ is defined by
$(\delta_s\otimes u)^*=\delta_{s\triangleleft u}\otimes u^i$.
\endproclaim
\demo{Proof} For the antipode:
$$
\align
ST(s\otimes \delta_{u})&=S(\delta_{u^{-1}}\otimes s^{-1})
=\delta_{(u^{-1}\triangleleft s^{-1})^{i}}\otimes (u^{-1}\triangleright s^{-1})^{i}\\
&=\delta_{((s\triangleright u)^{-1})^{i}}\otimes ((s\triangleleft u)^{-1})^{i}
\ ,\\
TS(s\otimes \delta_{u})&=T((s\triangleleft u)^{i}\otimes \delta_{(s\triangleright u)^{i}})
=\delta_{((s\triangleright u)^{i})^{-1}} \otimes ((s\triangleleft u)^{i})^{-1}\ .
\endalign
$$
 For the $*$-operation:
$$
\align
(*\circ T)(s\otimes \delta_{u})&=*(\delta_{u^{-1}}\otimes s^{-1})
=\delta_{u^{-1}\triangleleft s^{-1}}\otimes (s^{-1})^{i}
=\delta_{(s\triangleright u)^{-1}}\otimes (s^{-1})^{i}\ ,\\
(T\circ *)(s\otimes \delta_{u})&=T(s^{i}\otimes\delta_{s\triangleright u})
=\delta_{(s\triangleright u)^{-1}}\otimes (s^{i})^{-1}
=\delta_{(s\triangleright u)^{-1}}\otimes (s^{-1})^{i}\ .\square\endalign
$$
\enddemo

\proclaim{Proposition 8.6}
The maps $T$ and $T_J$ preserve the unit and counit.
\endproclaim
\demo{Proof} For the unit:
$$
\align
T\eta(j\otimes \delta_{n})&=T(\sum_{z\in G:zz^{i}=n} j\otimes \delta_{z})
=\sum_{z\in G:zz^{i}=n}  \delta_{z^{-1}}\otimes j^{-1}\ ,
\\
\eta T_J(j\otimes \delta_{n})&=\eta(\delta_{n^{-1}}\otimes j^{-1})
=\sum_{a\in M:aa^{i}=n^{-1}} \delta_{a}\otimes j^{-1}\ ,\endalign
$$
and these are the same by putting $a=z^{-1}$.
For the counit: If $u\in J_{G}$,
$$
\align
\epsilon T(s\otimes \delta_{u})&=\epsilon(\delta_{u^{-1}}\otimes s^{-1})
=\delta_{u^{-1}}\otimes s^{-1}(s^{-1})^{i}\ ,\\
&=\delta_{u^{-1}}\otimes (s^{i}s)^{-1}\ ,\\
T_J\epsilon(s\otimes \delta_{u})&=T_J(ss^{i}\otimes \delta_{u})\\
&=\delta_{u^{-1}}\otimes (ss^{i})^{-1} 
=\delta_{u^{-1}}\otimes (s^{i}s)^{-1}\ .\endalign
$$
If $u\notin J_{G}$ then both expressions
will give zero. \quad$\square$
\enddemo

\proclaim{Theorem 8.7}
The almost Hopf algebra $H=kM\subord k(G)$ is self dual by the map 
$H@>S >>H@>T >>H'$. 
\endproclaim
\demo{Proof} We have seen that both $S$ and $T$ reverse the order of 
the product and coproduct, and preserve the unit, counit and antipode.
Further both $S$ and $T$ are invertible.
\quad$\square$
\enddemo

\enddocument

\end{document}